\documentclass[12pt]{article}
\usepackage{latexsym}
\usepackage{amssymb}
 
\newcommand{\proof}{{\noindent \bf Proof. }}
\newtheorem{thm}{Theorem}
\newtheorem{defi}{Definition}

\newtheorem{prop}{Proposition}

\newcommand{\x}{{\bf x}}
\newcommand{\y}{{\bf y}}
\newcommand{\z}{{\bf z}}
\newcommand\mbf[1]{\mbox{\boldmath$#1$}}
\newcommand\msbf[1]{\mbox{\boldmath\scriptsize$#1$}}

\newcommand{\C}{{\cal C}}
\newcommand{\B}{{\cal B}}
\newcommand{\D}{{\cal D}}

\newcommand{\X}{{\cal X}}
\date{}

\begin{document}
\begin{titlepage}
\title{\bf ZERO-ERROR CAPACITY OF BINARY CHANNELS WITH MEMORY}
{\author{{\bf  G\'erard Cohen}
\\{\tt cohen@telecom-paristech.fr}
\\ENST
\\ FRANCE
\and{\bf Emanuela Fachini}
\\{\tt fachini@di.uniroma1.it}
\\''La Sapienza'' University of Rome
\\ ITALY
\and{\bf J\'anos K\"orner}
\thanks{Department of Computer Science, University of Rome, La Sapienza, 
via Salaria 113, 00198 Rome, ITALY}
\\{\tt korner@di.uniroma1.it}
\\''La Sapienza'' University of Rome
\\ ITALY}} 

\maketitle
\begin{abstract}

We begin a systematic study of the problem of the zero--error capacity of noisy binary channels with memory and solve some of the non--trivial cases.
\end{abstract}
\end{titlepage}

\section{Introduction}

Zero--error capacity problems in information theory represent an important challenge to the combinatorialist. They originate in the seminal paper of Claude Shannon \cite{Sh}. Beyond their 
relevance for Shannon Theory in itself, they offer a rather significant framework and proof techniques in asymptotic extremal combinatorics \cite{CK}, cf. also \cite{GKV}.
To our knowledge, the first paper about the zero-error capacity of channels with memory was written by Ahlswede, Cai and Zhang \cite{AZ}. They focus their attention on a rather particular long--time memory
channel model they call the enlightened dictator channel. A simpler but isolated model with long memory was considered in \cite{FK}. Our objective here is quite different. As a first step towards a systematic study of the zero-error capacity of channels with short memory, we are investigating the case in which the memory is of order 1 and the input  alphabet of the channel is binary. In other words, we consider the 4-element set of all the pairs of consecutive binary inputs to 
such a channel. The channel is then defined by specifying which pairs of these input sequences of length 2 can be distinguished at the channel output. It is easily seen that 
even if two channels have the same graph of confusability, their capacity can be different. This is true already when the graph has just one edge, as it was shown in \cite{AZ}.
Yet the problem of determining the zero-error capacity of channels whose confusability graph has a single edge is very easy and requires no new mathematical ideas. If, however, the graph has more edges, our problem becomes more challenging and we do not know all the answers. 

Formally, our problem is defined in terms of a finite set $\X$, the input alphabet of our channel and a confusability graph $G$ whose vertex set is the set $\X^2$ of consecutive pairs of 
{\em distinct} elements of the alphabet. The edge set $E(G)$ of the graph is an arbitrary but fixed subset of ${\X^2 \choose 2}.$ The sequences $\x \in \X^n$ and $\y \in \X^n$ are 
{\em distinguishable} for the channel graph $G$ if there is an $i\in [n-1]$ with $$\{x_ix_{i+1}, y_iy_{i+1}\}\in E(G).$$
(Restricting the edge set to pairs of distinct vertices is important for the intuition. Dropping this condition one encounters mathematically amusing problems 
of a different nature, cf. \cite{CFK}.)
Let $M(G, n)$ be the largest cardinality of a set $C\subseteq \X^n$ any two distinct sequences of which are distinguishable for the channel graph $G.$ We call 
$$C_M(G)=\limsup_{n \rightarrow \infty} {1 \over n} \log M(G, n)$$
the Shannon-Markov capacity of the graph $G.$ Except for the terminology, this model was introduced in \cite{AZ}. 

If $\X=\{0,1\}$ and the graph $G$ has just one edge, then 
\prop \noindent \rm{\cite{AZ}}

{\em $C_M(G)=1/2 $ if the two vertices of the edge of $G$ have Hamming distance 2 and $C_M(G)=\log{1+\sqrt{5} \over 2}$
 if the corresponding Hamming distance is 1.}
 
 The most immediate generalisation of this elementary result is obtained when we consider the graph $G$ to be the complete graph on 3 vertices. In the two 
 different cases of one--edge graphs one got different solutions precisely because of the different Hamming--distance of the vertices of the graph. We shall 
 see that for larger graphs such as the triangle we get different results even in case of isomorphic graphs with an isomorphism maintaining the pairwise Hamming distance of corresponding vertex pairs. (Such graphs are called isometric.)

\section{Triangles}

Let the graph $F$ be complete with vertex set $\{0,1\}^2-\{11\}.$ We have 

\begin{thm} \label{thm:11}
$$C_M(F)=\log \alpha\approx 0.878$$
where $\alpha^{-1}$ is the only positive root of  the equation 
$$x+x^2+x^3=1.$$
\end{thm}

\proof

We consider the set $\C_n^*=\{0, 01, 011\}^* \cap \{0,1\}^n$, for an arbitrary but fixed value of $n$. Note that since none of the strings $0, 01$ and $011$ is a postfix of an other one, every binary string in $\C_n^*$ has a unique decomposition into a sequence composed of disjoint substrings $0, 01$ and $011.$ By a well--known classical result of Shannon, 
(cf. e.g. Lemma 4.5 in \cite{CK}), 
$$\liminf_{n \rightarrow \infty} {1 \over n} \log |\C_n^*| \geq \alpha,$$
where $\alpha$ is the constant defined in the statement of this theorem.
We further claim that (a sufficiently large part of) the set $\C_n^*$ consists of pairwise distinguishable sequences for the channel graph $F.$ To see this, consider two distinct sequences, $\x\in C_n^*$ and $\y \in C_n^*.$ Because of the postfix--freeness of the set of mini--strings $ \{0, 01, 011\}$ these two strings have a different decomposition into mini--strings. Let $i \in [n]$ be the first coordinate of two different elements (ministrings) of $\{0, 01, 011\}$ appearing in the two different strings. Suppose first that one of these two ministrings is the singleton $0.$
Without restricting generality, we can suppose that this happens in the i'th coordinate of $\x.$ But then, since all the mini strings start with a zero, we have 
$x_ix_{i+1}=00.$ This implies that although $y_i=0$, $y_{i+1}=1$ which establishes our claim. Suppose next that none of the two different mini strings is $0.$
Then we have $x_ix_{i+1}=y_iy_{i+1}=01.$ Since these strings are the prefixes of two different ministrings, exactly one of $x_ix_{i+1}x_{i+2}$ and 
$y_iy_{i+1}y_{i+2}$ must equal 
$011.$ Suppose, without restricting generality, that it is $\x.$ This implies that $x_{i+2}x_{i+3}=10$ while $y_{i+2}=0$ (since it is the first digit of a new ministring) and this provides the desired 
difference, unless the ministring $011$ in question appears in the last three coordinates, i.e., $i+2=n.$ Hence this can only occur if all the other appearances of 
the ministring $011$ in the two strings $\x$ and $\y$ are in coinciding positions. This would imply that our two strings $\x$ and $\y$ have a different number of 
occurrences of the ministring $011.$ In order to exclude this let us partition $\C_n^*$ into classes according to the number of occurrences of $011$ in its strings.
The number of these classes is at most $\lfloor n/3 \rfloor.$ Let $\C_n$ be a class of maximum cardinality of our partition. Then
$$|\C_n|\geq {3|\C_n^*| \over n}$$
and $\C_n$ has all the properties we need. This proves
$$C_M(F)\geq \alpha.$$
In order to establish inequality in the opposite direction, let us consider the set $\D_n \subseteq \{0,1\}^n$ of all binary strings of length $n$ without three consecutive bits equal to 1. Clearly, $\C_n\subseteq \C_n^* \subseteq \D_n.$ We claim that 
$$|\D_n| \leq  3|C_n^*|.$$
As a matter of fact, $\D_n$ is the union of a set of strings $\D_n^{(1)}$ of strings beginning with a 1 with the set $\C_n^*.$ However, the strings in $\D_n^{(1)}$ can be obtained from strings in $C_{n-1}^*$ either by adding $1$ as a prefix to each of them or from strings in $C_{n-2}^*$ by adding a prefix $11.$ Hence 
$$\liminf_{n \rightarrow \infty} {1 \over n} \log |\D_n| \geq \alpha.$$
To conclude the proof, it is therefore enough to show that to every set of strings in $\{0,1\}^n$ being pairwise distinguishable for the channel $F$ there corresponds 
one of the same cardinality contained in $\D_n.$ To do so, consider an arbitrary set $\B\in \{0,1\}^n$ of pairwise distinguishable strings for our channel. Let 
$\x$ be an arbitrary string in $\B$ that contains at least three consecutive 1's. (If there is none, we are already done.) Let us replace an arbitrary substring of three 
consecutive 1's in $\x$ by the substring $101$ and let $\z$ be the string so obtained. Suppose that the middle coordinate of the three is $i \in [n].$
Clearly, $\z \notin \B,$ since it is not pairwise distinguishable from $\x.$
On the other hand, let $\y\not = \x$ be an arbitrary string from $\B.$ It is obvious that wherever there are two consecutive coordinates guaranteeing the pairwise dishinguishability of $\x$ and $\y,$ they will do it also for $\x$ and $\z$, since $\z$ differs from $\x$ only in the $i$'th coordinate and the latter appears only in the two 
2--length substrings $x_{i-1}x_i=x_ix_{i+1}=11.$ Thus replacing $\x$ by $\z$ in $\B$ leaves us with a good construction of the same cardinality as $\B.$ Iterating this 
procedure we eventually arrive at a subset of $\D_n$ as claimed.

\hfill$\Box$

Obviously, the problem has the same answer if the graph $F$ of the channel is complete with vertex set $\{0,1\}^2-\{00\}.$ To see this, it suffices to switch $0$ and 
$1$ in the previous theorem. Things change, however, for the complete graph $G$ with vertex set $\{0,1\}^2-\{10\}.$ Our next result shows that, somewhat 
surprisingly, the zero--error capacity of this channel is different. 
\begin{thm} \label{thm:10}

Let $G$ be the complete graph with vertex set $\{ 0,1\}^2-\{10\}$. 
We have
$$ C_M(G)=\log \beta \approx 0.849 $$
where $\beta^{-1}$ is the only positive root of the equation 
$$x+\frac {x^2}{1-x^2}=1.$$
\end{thm}

\proof

In order to obtain the claimed lower bound on $C_M(G)$ let $\C_n$ be the set of those sequences from $\{0,1\}^n$ in which every run of 1's has an odd length. A run is a 
maximal sequence of consecutive 1's. More precisely, it is a sequence of 1's which is not properly contained in a larger sequence of the same kind. Further, let each of the sequences in $\C_n$ have $0$ as their first coordinate. We claim that the strings of 
$\C_n$ are pairwise distinguishable for the channel graph $G.$ To see this, let $\x \in \C_n, \y \in \C_n$ be arbitrary but different. Let $j \in [n]$ be the first coordinate in which these two strings differ. Without restricting generality suppose that $x_j=0.$ Let further $i\leq j$ be the first coordinate of the run of 1's of $\y$ to which $y_j$ belongs. Suppose first that $j=i.$ In this case both $\x$ and $\y$ have a zero in the preceding coordinate and thus in these two coordinates the two sequences differ in the prescribed manner; we have $00$ in $\x$ and $01$ in $\y.$ Suppose next that $i<j.$ By our hypothesis 
$x_j$ is the first zero after a run of 1's in $\x$. Since all runs of 1's in our strings have an odd length, we conclude that also $y_{j+1}=1.$ But then in the coordinate pair 
$(j, j+1)$ our two strings differ in a pair of adjacent vertices of the graph $G$ as claimed. 

The set  $\C_n$ is the intersection of the sets $\{0,1\}^n$ and $\{0, 01, 0111, 011111,\dots\}^*.$ By Shannon's already cited classical theorem the cardinality of $\C_n$ satisfies
$$\lim_{n \rightarrow \infty}\frac{1}{n}\log |\C_n|=\beta$$
establishing the promised lower bound of the statement of our theorem. To explain this in somewhat more detail, by Shannon's theorem we know that for every fixed $k$ the cardinality of the set $\C_{n,k} \subseteq C_n$, defined as
$$\C_{n,k}=\{0,1\}^n \cap \{0, 01,\dots, 01^{2k+1}\}^*$$ satisfies
$$\lim_{n \rightarrow \infty}\frac{1}{n}\log |\C_{n,k}|=\log\beta_k$$
where $\beta^{-1}$ is the only positive root of the equation 
$$x+x^2+x^4+\dots +x^{2k+2}=1$$ and $1^{2k+1}$ denotes a binary string of length $2k+1$ containing no zero. Clearly, as $k$ goes to infinity, $\beta_k$ converges 
to $\beta.$

In order to prove the converse result, our upper bound, we denote by $\B_n$ the set of all those binary sequences of length $n$ every run of which has an odd length,
but now the first coordinate of a sequence might be 1. Hence, $\C_n \subset \B_n.$ In fact, it is easy to see that
those sequences in $\B_n$ which do not belong to $\C_n$ give rise to different sequences from $\C_{n+1}$ by adding a prefix 0 before their first coordinate.
Hence we have
$$|\B_n| < |\C_n|+|\C_{n+1}|< 3|\C_n|, $$
where the last inequality holds for $n$ sufficiently large considering that for such $n$ 
$$|\C_{n+1}|<2^{0.85}|\C_n|<2|\C_n|.$$

We now define a function $f_n: \{0,1\}^n\rightarrow \B_n$. Let $f_n(\x)=\x$ if $\x \in \B_n.$ For a sequence $\x \not \in \B_n$ let its image by $f_n$ be the 
sequence obtained from $\x$ by substituting the last 1 in every run of even length by a 0. Thus the image of every binary sequence of length $n$ is in 
$\B_n$ as claimed. Let us consider the partition of $\{0,1\}^n$ generated by the function $f_n.$ It is clear that if a pair of sequences $\{\x, \y\}$ is in the same 
class of this partition, i.e., if $f_n(\x)=f_n(\y),$ then the two sequences  $\x$ and $\y$ do not satisfy our condition to be distinguishable for the channel represented by $G.$ (This is easy to see. The two sequences have the first coordinates of their respective runs  of 1's in the same places. Their corresponding runs, those beginning in the same coordinate, have lengths differing by at most one. The not coinciding last coordinates of two corresponding runs are therefore consecutive and produce a difference where one of the sequences has 11 and the other one has 10.) Let now $\D\subset \{0,1\}^n $ be an optimal code of length $n$ for $G.$ Thus
$$|\D|=M(G,n).$$
By the foregoing, the function $f_n$ is injective on $\D.$ Hence, 
$$|M(G,n)|=|\D|\leq |\B_n|<3|\C_n|.$$
This, using the asymptotics of $\C_n$ from the first part of our proof, establishes the claimed upper bound and thus completes the proof.

\hfill$\Box$

We have seen that $C_M(G)<C_M(F)$. The proofs for the two capacities are different and, unfortunately, there seems to be a lack of general methods to tackle these apparently simple problems. 
Just to go one step further in this exploration, let us consider the two non--isomorphic cases associated with the complete bipartite graph $K_{1,3}.$ 
\begin{thm} \label{thm:star1}

Let $L$ be the 
graph with vertex set $\{ 0,1\}^2$ all of whose three edges are incident to $(0,0).$ We claim
 
$$ C_M(L)= \log \gamma \approx 0.81 $$
where $\gamma^{-1}$ is the only positive root of the equation 
$$x+\frac {x^3}{1-x}=1.$$
\end{thm}

\proof

Let $g: \{ 0,1\}^2 \rightarrow \{ 0,1\}$ be 
$g(0,0)=0$ and have the value $1$ for the remaining three binary pairs. 
We define a function $f_n: \{ 0,1\}^n \rightarrow \{ 0,1\}^n$ as follows. Let $f_n(\x)$ have its first coordinate equal to that of $\x$, and for $i>1$ its 
$i-$th coordinate be equal to $g(x_{i-1}, x_i)$. Let us first restrict attention to the subset $C_n$ of the domain of $f_n$ which contains 
strings with $0$ in their first coordinate. Then it is clear that although that the restriction of the function $f$ to $C_n$ is not injective, we have $f_n(\x)\not=f_n(\y)$ precisely when the strings in  $\x$ and $\y$ of $C_n$ are distinguishable for the channel graph $L.$ Let us denote by $D_n$ the true range of $f_n.$ Next we can partition $C_n$ into at most $n$ classes such that a string belongs to the $j-$th class if its first digit equal to $0$ appears in its $j-$th coordinate. This shows that the largest cardinality of a good code for our channel is between $|D_n|$ and $n|D_n|$. Note that $D_n$ consists of precisely those strings from $ \{ 0,1\}^n$ which have first coordinate $0$ and do not have isolated $1'$s. In view of Shannon's already cited theorem, the cardinality of the set $ D_n$ satisfies
$$\lim_{n \rightarrow \infty} \frac{1}{n}\log |D_n|=\log \gamma.$$

\hfill$\Box$

Finally, let $Q$ be once again the complete bipartite graph $K_{1,3}$ but let this time $(0,1)$ be the vertex of degree $3.$ Then 

\begin{thm} \label{thm:star2}

Let $Q$ be the 
graph with vertex set $\{ 0,1\}^2$ all of whose three edges are incident to $(0,1).$ We claim
 
$$ C_M(Q)= \log \frac{1+\sqrt{5}}{2}.$$

\end{thm}

\proof
As a lower bound on $C_M(Q)$, the statement follows from Proposition 1, since the capacity of $Q$ is lower bounded by that of its subgraph having a single edge 
with endpoints $(0,1)$ and $(0,0).$ To establish a matching upper bound, we proceed as in the previous theorem. We define the function 
$g: \{ 0,1\}^2 \rightarrow \{ 0,1\}$ by setting $g(0,1)=1$ and having the value $0$ for the remaining binary pairs. We next define a function 
$f_n: \{ 0,1\}^n \rightarrow \{ 0,1\}^n$ by letting $f_n(\x)$ have its first coordinate equal to that of $\x$, and for $i>1$ its 
$i$-th coordinate be equal to $g(x_{i-1}, x_i)$. It should be clear that this function is injective on the codewords of an $n$--length block code for $Q.$ 
It is equally clear that the values of the function are binary strings without consecutive $1'$s, the so--called Fibonacci sequences, which completes the proof.

\hfill$\Box$

To widen our horizon we will show how the previous questions can be regarded as capacity problems for memoryless channels with an input constraint. 

\section{Input constraints}

The concept of capacity of memoryless channels with a constrained input arises naturally when one deals with the compound channel, i.e., a channel whose unknown 
transmission probability matrix belongs to a finite set of possible alternatives. More precisely, in this case the key ingredient in the formula for the zero-error capacity of the channel (with an informed decoder) is a concept introduced by
Csisz\'ar and K\"orner \cite{CKW}.  They needed the notion of zero--error channel capacity  ``within a fixed type''. This is based on the notion of {\em types}.

\begin{defi}
The  type of a sequence $\mbf x\in V^n$ is the probability distribution
$P_{\msbf x}$ on $V$ defined
by $$P_{\msbf x}(a)={{|\{i: x_i=a\}|}\over n}, \mbox{ for all } a\in V.$$
For a fixed distribution $P$ on $V$ and $\varepsilon>0$, we say that $\mbf x\in
V^n$ is {\em $(P,\varepsilon)$-typical} if, for all $a\in V$, we have 
$|P_{\msbf x}(a)-P(a)|< \varepsilon.$ 
\end{defi}

Capacity is the asymptotic speed of growth of the largest clique in the powers of the graph $G$. 
Let $G=(V(G), E(G))$ be a simple graph. Thus $E(G)\subseteq {V(G)\choose 2}.$ The graph $G^n$ has as 
vertices the sequences of length $n$ of the vertices of $G.$ We have 
$$\{\x, \y\} \in E(G^n) \quad \hbox{if}\quad \exists i\in [n] \quad \hbox{with}\quad \{x_i,y_i\}\in E(G).$$
The cardinality of the largest complete subgraph in a graph $G$ is denoted by $\omega(G).$
\begin{defi}{\rm (cf. \cite{CKW})}
The (logarithmic) Shannon capacity within type $P$ of a (finite) graph $G$ with vertex set $V$ is 
$$C(G,P)=\lim_{\varepsilon\to 0}\limsup_{n\to\infty}{1\over
  n}\log\omega(G^n(P,\varepsilon)),$$
where $G^n(P,\varepsilon)$ denotes the graph induced by $G^n$ on the
$(P,\varepsilon)$-typical sequences in $V^n$. 
\end{defi}

What we want to consider here is (a special case of) the extension of the previous definition to (topological) Markov types. More precisely, let $P$ be a directed graph with vertex set $V$. The edge set of $P$ is an arbitrary subset of $V^2.$ In particular, loops are not excluded. We denote by $V^n(P)$ the set of those sequences $\x\in V^n$ for which 
$$(x_i, x_{i+1})\in E(P) \quad \hbox{for every} \quad i<n.$$ We denote by $G^{n,P}$ the graph $G^n$ induces on $V^n(P).$
This graph will play the role of a type in our present context.
We introduce
\begin{defi}
$$C_P(G)=\limsup_{n\rightarrow \infty}\frac{1}{n} \log \omega(G^{n,P})$$
and call it the (logarithmic) zero-error capacity of the channel $G$ within the topological Markov type $P.$
\end{defi}

It should be clear that this generalises our previous concept of zero--error capacity, the Shannon--Markov capacity of a graph. In order to explain this, we will 
show how $C_M(F)$ can be redefined in this setting. We set $V=\{0,1\}^2$ and define an edge in $P$ pointing from $(x_1,x_2)$ to 
$(y_1,y_2)$ if $x_2=y_1.$  With this definition every string in $\{0,1\}^n$ gives rise to a string of length $n-1$  of vertices from $V(P)$ in a bijective manner.
With this definition, we have the equality
$$M(F,n)=\omega(G^{n-1,P}).$$

It is interesting to extend the  set--up of channel codes within a fixed topological Markov type from simple graphs to directed graphs since this will 
allow us to integrate into the topic of capacity of graphs (and graph families) some previously scattered and apparently unrelated problems from extremal 
combinatorics.

A very natural generalisation of Shannon's graph capacity for directed graphs was introduced in \cite{GKVfirst} by the name Sperner capacity. This concept was the key to the solution of 
a well--known open problem of R\'enyi on the largest family of pairwise qualitatively independent $k$--partitions of an $n$--set, cf. \cite{GKV} and several other problems 
in an outside information theory, cf. Chapter 11 of the book \cite{CK}. The following definitions are from \cite{GKVfirst}.

Let $G$ be a directed graph with vertex set $V=V(G).$ A set $C\subseteq V(G)$ is said to induce a symmetric clique in $G$ if every ordered pair of distinct vertices 
from $C$ is an edge in $G.$ Let us denote by $\omega_s(G)$ the largest size of a symmetric clique in $G.$ Next we define the power graphs of a directed graph.
For any natural number $n$ the graph $G^n$ has vertex set $V(G^n)=[V(G)]^n$. There is an edge from $\x \in V(G^n)$ to $\y \in V(G^n)$ if at least in one coordinate, 
$i \in [n]$ we have $(x_i, y_i) \in E(G)$, just like in the undirected case. We define
\begin{defi}\noindent\rm{\cite{GKVfirst}}

The (logarithmic) Sperner capacity of the digraph $G$ is the always existing limit
$$Sp(G)=\lim_{n\rightarrow \infty}=\frac{1}{n}\log \omega_s(G^n).$$
\end{defi}

This name is justified by the fact that if $S$ is the single edge graph on a two-element vertex set, then $\omega_s(G^n)$ is the largest cardinality of a family of 
subsets of $[n]$ such that none of the member sets contains an other one. This observation shows that Sperner's theorem \cite{Sp} is strongly related to 
the problem area around Shannon's zero--error capacity. It is straightforward to extend the concept of capacity of a graph in fixed topological Markov chain $P$ to directed graphs by considering Sperner capacity. More precisely, given the digraphs $G$ and $P$ with the same vertex set $V$ we set
$$Sp(G,P)=\limsup_{n\rightarrow \infty}=\frac{1}{n}\log \omega_s(G^n,P)$$
where $\omega_s(G^n,P)$ is the largest cardinality of a symmetric clique $G^n$ induces on the set $V^n(P).$

In \cite{CFK} we have introduced the following very elementary problem. Let $\sf{F}_n$ be the 
set of all the binary sequences of length $n$ without 1's in consecutive positions. (The cardinality of this set is the classical example for the standard Fibonacci sequence.)
Consider these binary sequences as the characteristic vectors of subsets of the set $[n]$ in the usual manner. We ask for the maximum cardinality of a Sperner family they contain. The hitherto sharpest result on this problem is due to Victor Falgas-Ravry \cite{FR}. It is immediate to realize that this problem has a natural formulation 
in our set--up. To set ideas, let $Fib(n)$ be this largest cardinality. Consider the directed graph $P$ with vertex set $\{0,1\}$ and edge set 
$\{(0,0), (0,1), (1,0)\}.$ Let further $G$ be a directed graph with vertex set $\{0,1\}$ and the single edge $(0,1).$ With this notation we have
$$Fib(n)=\omega_s(G^n,P).$$
Also, it is trivial that
$$Sp(G,P)=\log\frac{1+\sqrt{5}}{2},$$
with$Fib(n)$ having the same exponential asymptotics.

Open problems abound. We conclude by just one. For an arbitrary natural number $k$ let $K_k$ be the symmetric clique with no loop edges. We are interested in determining $Sp(G, K_k)$ for an arbitrary directed graph $G$ on the vertex set of $K_k$. As a matter of fact, this problem is interesting also in the case of a simple 
graph, and needs no new definition, since the Shannon capacity of a simple graph is equivalent to the Sperner capacity of the digraph obtained from it by replacing each of its edges by two directed edges between the same vertices, going in opposite directions. We believe that for the pentagon $C_5$
$$Sp(C_5, K_5)=1.$$

\end{document}